\documentclass[10pt,a4paper]{article}
\usepackage{amssymb}
\usepackage{amsmath}
\usepackage{amsfonts}
\usepackage{dsfont}
\usepackage{hyperref}
\usepackage[left=3cm,top=3cm,right=2cm,bottom=2cm]{geometry}

\begin{document}
\newtheorem{defi}{Definition}
\newtheorem{lema}{Lemma}
\newtheorem{coro}{Corollary}
\newtheorem{teor}{Theorem}
\newtheorem{prop}{Proposition}

\title{The Ring of Integers in the Canonical Structures of the Plane}
\author{Jos\'e C. Cifuentes$^a$, Jo\~ao E. Strapasson$^b$, Ana C. Corr\^ea$^a$ \and Patr\'\i cia M. Kitani$^a$ \\ \small $^a$ Department of Mathematics Federal University of Paran\'a \and \small $^b$ Institute of Mathematics, State University of Campinas }
\maketitle

\begin{abstract}
The \emph{canonical structures of the plane} are those that result, up to isomorphism, from the rings that have the form $\mathds{R}[x]/(ax^2+bx+c)$ with $a\neq 0$.That ring is isomorphic to $\mathds{R}[\theta]$, where $\theta$  is the equivalence class of x, which satisfies $\theta^2 = \left( -\dfrac{c}{a} \right) + \theta \left(-\dfrac{b}{a}\right)$. On the other hand, it is known that, up to isomorphism, there are only three canonical structures: the corresponding to $\theta^2 = -1$ (the complex numbers), $\theta^2 = 1$ (the perplex  or hyperbolic numbers) and $\theta^2 = 0$ (the parabolic numbers). This article copes with the algebraic structure of the rings of integers $\mathds{Z}[\theta]$ in the perplex  and parabolic cases by \emph{analogy} to the complex cases: the ring of Gaussian integers. For those rings a \emph{division algorithm} is proved and it is obtained, as a consequence, the characterization of the prime and irreducible elements.  
\end{abstract}  

\section{The Plane Canonical Structures}
The Cartesian plane $\mathds{R}^2$ supports a very rich family of algebraic structures,and one of the most important the complex numbers $\mathds{C}$. Starting from the vector sum on $\mathds{R}^2$, we may ask: what products may be defined in a way that is compatible with the sum? In analogy to the complex numbers, we may think the elements of $\mathds{R}^2$ as $z = x + \theta y$ with $x, y \in \mathds{R}$ and $\theta$ a new object such that $\theta^2 = \alpha +\theta\beta$, where $\alpha$ and $\beta$ are real constants. In that case, the product, defined distributively with respect to the sum, has the following form: $$(x_1 + \theta y_1)(x_2 + \theta y_2) = (x_1 x_2 + \alpha y_1y_2) + \theta(x_1y_2 + x_2y_1 + \beta y_1y_2).$$

Despite to the infinity of possible values for $\alpha$ and $\beta$  it can be demonstrated, through the discriminant $D = \beta^2 + 4\alpha$, that there are, up to isomorphism, only three structures for $\mathds{R}^2$ which correspond to the values of $\theta^2 = -1$ (the elliptic case $D < 0$ ), $\theta^2 = 1$ (the hyperbolic case $D > 0$ ) and $\theta^2 = 0$ (the parabolic case  $D = 0$) \cite{4}. 
The mentioned discriminant results from analyzing the norm $\eta(x + \theta y) = \left(x + \dfrac{1}{2}\beta y\right)^2 - \dfrac{1}{4}Dy^2$, which is obtained from the  minimal polynomial of an element $z = x + \theta y$. That polynomial has the form $P(z) = z^2 - (2x + \beta y)z + (x^2 + \beta xy - \alpha y^2)$ and so, the \emph{trace} $\tau(z) = 2x + \beta y$ and the \emph{norm} $\eta(z) = x^2 + \beta xy - \alpha y^2 = \left(x + \dfrac{1}{2}\beta y\right)^2 - \dfrac{1}{4}Dy^2$. The elliptic case corresponds to the structure of the field of \emph{complex numbers} $\mathds{C}$ and its imaginary unit is denoted by $\theta = i$, the hyperbolic case corresponds to the ring of \emph{perplex  numbers} \cite{1} or \emph{hyperbolic numbers} $\mathds{H}$ and its imaginary unit will be denoted by $\theta = j$, and the parabolic case, not yet properly studied, corresponds to the ring $\mathds{P}$ of \emph{parabolic numbers} and whose imaginary unit is denoted by $\theta = k$. Any one of them can be denoted by $\mathds{R}[\theta]$. Only in the case of the complex number that we have a field; in the other cases the ring is not an integral domain, even though, it is a commutative ring with unit 1, and $\mathds{R}$ can be embedded in $\mathds{R}[\theta]$ in the usual way. On the other hand, it can be prove that, in the case of the ring $\mathds{P}$, the \emph{lexicographical order} of $\mathds{R}^2$ is consistent with the algebraic structure in the sense that we have the structure of an ordered ring, which makes $\mathds{P}$ an extension of $\mathds{R}$ that allows the existence of infinitesimals. In fact, all elements of the form ky, with real y,  are infinitesimals of $\mathds{P}$ \cite{6}.  

The aim of this paper is to analyze the structure of the integers ring $\mathds{Z}[\theta]$ for the hyperbolic and parabolic cases, that is, the rings of hyperbolic integers $\mathds{Z}[j]$ and of parabolic integers $\mathds{Z}[k]$, emphasizing the characterization of prime and irreducible elements.

This article's main purpose is to analyze the structure of the ring of integers $\mathds{Z}[\theta]$ for the hyperbolic and parabolic cases, that is, the rings of \emph{integers hyperbolic} $\mathds{Z}[j]$ and of \emph{integers parabolic} $\mathds{Z}[k]$, emphasizing the characterization of elements prime and irreducible.

This study begins with the proof of an ``appropriate" \emph{division algorithm} for those rings, allowing as particular case the ring of \emph{Gaussian integers} $\mathds{Z}[i]$. In the other cases, it is given a proper place to the zero divisors. The difficulties of translation and adaptation of the properties of $\mathds{Z}[i]$ in the other cases can be noticed right at the start, since the new rings are not, as we mentioned, integral domains and, as far as we know, there  exists no general theory of rings with an algorithm for division where the zero divisors play an essential role \cite{5}.  

Next we describe very briefly the structure of $\mathds{R}[\theta]$ in a unified form for all the three cases.  

In $\mathds{R}[\theta]$ it is possible to define, in analogy to the complex case, a \emph{conjugate element} and a \emph{norm} in the following way: for $z = x + \theta y$,  the conjugated of $z$ is given by $\bar{z} = x - \theta y$ and the norm of $z$ by $\eta(z) = x^2 - \theta^2y^2$; thus $z\bar{z} = \eta(z)$. Observe that, in the complex case $\eta(z) = x^2 + y^2 = |z|^2$, in the perplex or hyperbolic cases $\eta(z) = x^2 - y^2$, and in the parabolic case $\eta(z) = x^2$

It is important to point out that the concept of ``norm" adopted  here a generalization that we consider quite suitable for rings with zero divisors, namely: if $\mathcal{A}$ is a ring and $\mathcal{D}$ is the set of zero divisors  of $\mathcal{A}$ (including the zero of the ring), then, a \emph{norm} in $\mathcal{A}$ is a function $\eta : \mathcal{A} \longrightarrow \mathds{Z}$ such that (a) $\eta(a) = 0 \Longleftrightarrow a \in \mathcal{D}$; and (b) $\eta(ab) = \eta(a)\eta(b)$. The norm $\eta$ is said \emph{positive} if for all $a, \eta(a) \geq 0$.

If $\eta$ is a norm in $\mathcal{A}$, then, $\eta^+(a) = |\eta(a)|$, for all $a$ is said a positive norm in $\mathcal{A}$.

By means of these concepts we can express several properties of the algebraic structure of $\mathds{R}[\theta]$. Thus:  
\begin{enumerate}
\item $\overline{(z + w)} = \bar{z} + \bar{w}$ and $\overline{(zw)} = \bar{z}\bar{w}$.  
\item $\overline{\overline{z}} = z$, and $\bar{z} = z \Longleftrightarrow z \in \mathds{R}$.  
\item $\eta(z) = z\bar{z}$ and $\eta(zw) = \eta(z)\eta(w)$.

In particular, the last property expresses, in the case of integer values, that the sum of squares, difference of squares and perfect squares, are of the same type.  
\item \emph{Law of the Parallelogram}: $\eta(z + w) + \eta(z - w) = 2(\eta(z) + \eta(w))$.  
\item $z$ is \emph{invertible} $\Longleftrightarrow  \eta(z) \neq 0$ and, in that case, $z^{-1} = \dfrac{\bar{z}}{\eta(z)}$. 
\item $z$ is \emph{zero divisor} $\Longleftrightarrow \eta(z) = 0$. 

Denoting $\mathcal{D}$ as the set of divisors of zero of $\mathds{R}[\theta]$ we have:  
\begin{itemize}
\item In the case $\mathds{C}, \mathcal{D} = \{0\}$, that is, the origin of $\mathds{R}^2$.  
\item In the case $\mathds{H}, \mathcal{D} = \{x \pm \theta x | x \in \mathds{R}\}$, that is, the principal and secondary diagonals of the plane.  
\item In the case $\mathds{P}, \mathcal{D} = \{\theta y | y \in \mathds{R}\} (= \{z \in \mathds{R}[\theta] | z$ is infinitesimal$\}$), that is, the $y$ axis.
\end{itemize}

It can also be prove that the norm comes from an (indefinite) inner product given by:   
$$\langle z, w \rangle = x_1 x_2 + \theta^2 y_1y_2,$$
for $z = x_1 + \theta y_1$ and $w = x_2 + \theta y_2$.  
In that case we have:  
\item $\langle z, w \rangle = Re (z\bar{w})$, where $Re (z) = \dfrac{z + \bar z}{2}$, and $\langle z, z \rangle = \eta(z)$.  
\item \emph{Law of Polarization}: $\langle z, w \rangle = \dfrac{1}{4}(\eta(z + w) - \eta(z - w))$.  
\item \emph{Law of Cosines}: $\eta(z - w) = \eta(z) + \eta(w) - 2\langle z, w\rangle$.  
\item \emph{Inequality of Schwarz}:  \\
\begin{itemize}
\item In the case $\mathds{C}: \langle z, w\rangle^2 \leq \eta(z)\eta(w);$ 
\item In the case $\mathds{H}: \langle z, w\rangle^2 \geq \eta(z)\eta(w);$  
\item In the case $\mathds{P}: \langle z, w\rangle^2 = \eta(z)\eta(w)$.
\end{itemize}
\end{enumerate}

Finally, we have the following algebraic representation: $\mathds{R}[\theta] \cong \mathds{R}[x]/(x^2 - \theta^2)$, thus, $\mathds{C} \cong \mathds{R}[x]/(x^2 + 1)$, $\mathds{H} \cong \mathds{R}[x]/(x^2 - 1)$ and $\mathds{P} \cong \mathds{R}[x]/(x^2)$, which are particulars cases of the ring $\mathds{R}[x]/( ax^2 + bx + c)$ with a $\neq 0$. It is prove that if $D = b^2 - 4ac$, then, that last ring is isomorphic to $\mathds{R}[\theta]$, where $\theta$ is the equivalence class of x and $\theta^2 = -1$ if $D < 0, \theta^2 = 1$ if $D > 0$ and $\theta^2 = 0$ if $D = 0$.

The perplex  numbers, although they do not form a field, they do have a close similarity with to the complex numbers. Perplex numbers are related to the hyperbolic functions in the same way that the complex numbers are related to the circular (trigonometrics) functions. For example, it can be proved that all the perplex  numbers z that are not zero divisor admit a hyperbolic representation satisfying an analog of the \emph{Moivre's theorem}. Thus, for example, if $z = x + jy$ with $\eta(z) > 0$ and $x > 0$, there exists $\alpha \in \mathds{R}$  such that $z =\sqrt{\eta(z)} (\cosh \alpha + j \sinh \alpha)$, and if $n \in \mathds{Z}$ we have that $z^n = (\sqrt{\eta(z) })^n (\cosh n\alpha+ j \sinh n\alpha )$. Besides it ca be defined the \emph{ perplex exponential } function in the following way: $\exp z = e^x (\cosh y + j \sinh y)$, where the following \emph{ Euler formulas} are satisfied: $\cosh x = \dfrac{e^{jx} + e^{-jx}}{2}$ and $\sinh x = \dfrac{e^{jx} - e^{-jx}}{2j}$, which is a perplex  reformulation of the well known $\cosh x = \dfrac{e^x + e^{-x}}{2}$ and $\sinh x = \dfrac{e^x - e^{-x}}{2}$. In the parabolic case, such an analogy entails us to define the parabolic functions, \emph{the parabolic cosine} and \emph{parabolic sine}, in the following way: $\rm{cos} p\ x = 1$ and $\rm{sin} p\ x = x$ for all $x$ \cite{2} and \cite{3}.

\section{A Division Algorithm for $\mathds{Z}[\theta]$}  
From now on, we define the following positive norm, in an unified form, for the complex, perplex  and parabolic cases: for $z = x + \theta y \in \mathds{Z}[\theta], \eta^+(z) = |x^2 - \theta^2 y^2|$. We also denote $\mathcal{D}$ as the set $\mathcal{D} \cap \mathds{Z}[\theta]$, and as usual, $(z)$ will denote the principal ideal generated by $z \in \mathds{Z}[\theta].$
  
Regarding the Gaussian integers, as we already observed, $\mathcal{D} = (0)$. If $\theta = j$, then, $\mathcal{D} = \mathcal{D}^+ \cup \mathcal{D}^-$ where $\mathcal{D}^+ = \{x + jx | x \in \mathds{Z}\}$ and $\mathcal{D}^- = \{x - jx | x \in \mathds{Z}\}$, so, $\mathcal{D}^+$ and $\mathcal{D}^-$ are respectively the principal and secondary diagonals of $\mathds{Z}\times\mathds{Z}$. If $\theta = k$, then, $\mathcal{D} = \mathcal{D}_0 = \{ky | y \in \mathds{Z}\}$, and so, $\mathcal{D}_0$ is the $y$ axis of $\mathds{Z}\times\mathds{Z}$. 
 
Next, we analyze the structure of $\mathcal{D}$ and its relationship to the ideals of $\mathds{Z}[\theta]$.   
  
\begin{prop} $\mathcal{D}^+$, $\mathcal{D}^-$ and $\mathcal{D}_0$ are principal and prime ideals of the respective $\mathds{Z}[\theta]$.
\end{prop}  
{\bf Proof:} It is trivial the fact that they are ideals of $\mathds{Z}[\theta]$. It can also be easily proved that 
$\mathcal{D}^+ = (1 + j), \mathcal{D}^- = (1 - j)$ and $\mathcal{D}_0 = (k)$. 

Now we have to prove that $\mathcal{D}^+$ is prime: suppose that $zw \in \mathcal{D}^+$ with $z = x + jy$ and $w = r + js$, so, $(xr + ys) + j(xs + yr) \in \mathcal{D}^+$. Therefore, $xr + ys = xs + yr$; thus $(x - y)(r - s) = 0$; therefore, $x = y$ or $r = s$, so, $z \in \mathcal{D}^+$ or $w \in \mathcal{D}^-.$   \hfill $\Box$
  
\begin{prop} If $\mathcal{I}$ is an ideal of $\mathds{Z}[j]$ with $\mathcal{I} \subset \mathcal{D}$, then, $\mathcal{I} \subset {D}^+ $or $\mathcal{I} \subset \mathds{D}^{-}.$  
\end{prop}  
{\bf Proof:} Suppose that $\mathcal{I} \not\subset \mathcal{D}^+$ and $\mathcal{I} \not\subset \mathcal{D}^-$. Then, there are $a , b \in \mathcal{I}$ such that $a \not\in \mathcal{D}^+$ and $b \not\in \mathcal{D}^-$, in particular $a \neq 0 $ e $b \neq 0$. Since $\mathcal{I}\subset \mathcal{D}$, we have $a \in \mathcal{D}^-$ and $b \in \mathcal{D}^+$, that is, $a = x - jx$ and $b = y + jy$; thus, $a + b = (x + y) - j(x - y) \in \mathcal{I} \subset \mathcal{D}$; therefore, $x + y = - (x - y)$ or $x + y = x - y$, from which, $x = 0$ or $y = 0$, that is, $a = 0$ or $b = 0$, a contradiction.   \hfill $\Box$

\begin{teor}(Division Algorithm in $\mathds{Z}[\theta]$)
Let $a, b \in \mathds{Z}[\theta]$ with $\eta^+(b) \neq 0$. Then, there exist $\gamma, \rho \in \mathds{Z}[\theta]$ such that $a = \gamma b + \rho$ with $\eta^+(\rho) < \eta^+(b)$.  
\end{teor}
  
\noindent{\bf Proof:} We need to find $\gamma \in \mathds{Z}[\theta]$ such that $\eta^+(a - \gamma b) < \eta^+(b).$
Since $\eta^+(b)\eta^+\left(\dfrac{a}{b} - \gamma\right) = \eta^+(a - \gamma b)$, we need to find $\gamma \in \mathds{Z}[\theta]$ such that $\eta^+\left(\dfrac{a}{b} - \gamma\right) < 1$. We have that $\dfrac{a}{b} = x + \theta y$ with $x, y \in \mathds{Q}$.
   
Let $r, s \in \mathds{Z}$ such that $|x - r| \leq \dfrac{1}{2}$ and $|y - s| \leq \dfrac{1}{2}$. Suppose that $\gamma = r + \theta s$ and $\rho = a - \gamma b$. Hence $a = \gamma b + \rho$  and $\eta^+(\rho) = \eta^+(b)\eta^+\left(\dfrac{a}{b} - \gamma\right) = \eta^+(b)\eta^+((x - r) + \theta (y - s)) = \eta^+(b)|(x - r)^2 - \theta^2(y - s)^2| \leq \eta^+(b) ((x - r)^2 + (y - s)^2) \leq \eta^+(b)\left(\dfrac{1}{4} + \dfrac{1}{4}\right) = \eta^+(b)\left(\dfrac{1}{2}\right)< \eta^+(b).$
  
The following proposition is the analogous to rings of principal ideals, in the case of integral domains.    \hfill $\Box$ 

\begin{prop}
If $\mathcal{I}$ is an ideal of $\mathds{Z}[\theta]$ with $\mathcal{I} \not\subset \mathcal{D}$, then, there exists $\alpha \in \mathds{Z}[\theta]$ such that $\eta^+(\alpha)\neq 0$ and $\mathcal{I} = (\alpha) + \mathcal{I}\cap \mathcal{D}$. \label{prop 3}
\end{prop}  

\noindent {\bf Proof:} Since $\mathcal{I} \not\subset \mathcal{D}$, there exists $\alpha \in \mathcal{I}$ such that $\eta^+(\alpha)\neq 0$. Let $C = \{\eta^+(z) | z \in \mathcal{I}\mbox{ and }z \notin \mathcal {D}\}, m = {\rm min} C$ and $\alpha \in \mathcal{I}$ such that $\eta^+(\alpha) = m$.  

We will prove that $\mathcal{I} = (\alpha) + \mathcal{I} \cap \mathcal{D}$.  

Let $z \in \mathcal{I}$. Since $\eta^+(\alpha)\neq 0$,  there exist $\gamma, \rho \in \mathds{Z}[\theta]$ such that $z = \gamma\alpha + \rho$ with $\eta^+(\rho)<\eta^+(\alpha) = m.$  

Since $\rho= z -\gamma \alpha \in \mathcal{I}$, then, $\eta^+(\rho) = 0$, because $m$ is minimum. So, $\rho \in \mathcal{D}$; therefore, $\rho \in \mathcal{I}\cap \mathcal{D}$ and $z \in (\alpha) + \mathcal{I}\cap\mathcal{D}$.  

If $z \in (\alpha) + \mathcal{I}\cap \mathcal{D}$, then, $z = \gamma\alpha +\rho$ for some $\gamma$ and $\rho$. Therefore, since $\alpha, \rho \in \mathcal{I}$, we have that $z \in \mathcal{I}$. \hfill $\Box$ 
  
Note that, as it is well known, $\mathds{Z}[i]$ is a ring of principal ideals. We observe as well, that the demonstration of  proposition \ref{prop 3} suggests, in the case of $\mathds{Z}[j]$, a modification the form of the ideal $\mathcal{I}$ in the following way: $\mathcal{I} = (\alpha) + \mathcal{I}\cap\mathcal{D}^+ + \mathcal{I} \cap \mathcal{D}^-$, where $\mathcal{I} \cap \mathcal{D}^+$ and $\mathcal{I}\cap \mathcal{D}^-$  are principal ideals of $\mathcal{D}^+$ and $\mathcal{D}^-$ respectively.

\section{Some Results about Units and Associated Elements}
One of the first results of this research, along with the identification of the algorithm of the division, is the characterization of the \emph{unit elements} of the ring $\mathds{Z}[\theta]$. We are going to see that the complex and perplex cases  are similar, although, in the parabolic case there is an essential difference.  

\begin{prop}
Let $a \in \mathds{Z}[\theta]$. The following statements are equivalent:  
\begin{description}
\item[i] $a$ is unit in $\mathds{Z}[\theta]$, that is,  $a \in U(\mathds{Z}[\theta])$.  
\item[ii] $\eta^+(a) = 1$.  
\item[iii] $a\in \{-1, 1, -\theta , \theta\}$ if $\theta$ is equal to $i$ or $j$, and $a\in  \{\pm 1 + \theta y | y \in \mathds{Z}\}$ if $\theta = k.$  
\end{description}
\end{prop}

\noindent {\bf Proof:}\\
($i \longrightarrow ii$): If $a$ is unit, then, there is $b$ such that $ab = 1$. Thus, $\eta^+(a) \eta^+(b) = \eta^+(ab) = \eta^+(1) = 1$, therefore, $\eta^+(a) = 1$.\\  
($ii \longrightarrow iii$): Suppose that $\eta^+(a) = 1$ and $a = x + \theta y$, then, $|x^2 - \theta^2y^2| = 1$. If $\theta^2 = -1$ then $x^2 + y^2 = 1$, where the solutions are $x = 0$ and $y =\pm 1$ or $y = 0$ and $x =\pm 1$, that is, $a \in \{-1 , 1 , -\theta , \theta\}$. If $\theta^2 = 1$ then $x^2 - y^2  = 1$ or $x^2 - y^2 = -1$, where the solutions are also $x = 0$ and $y =\pm 1$ or $y = 0$ and $x =\pm 1$, that is, $a \in \{-1 , 1 , -\theta , \theta\}$. If $\theta^2 = 0$ then $x^2 = 1$, that is, $x = \pm 1$ and $y$ is any value; therefore, $a \in \{(\pm 1 + \theta y | y \in \mathds{Z}\}$.\\  
($iii \longrightarrow  i$): In the cases $\theta = i$ or $\theta = j$ the elements $1, -1, \theta$ and $-\theta$ are units. In the case $\theta = k$, the elements $\pm 1 + \theta y$ have $\pm 1 - \theta y$ as inverse.   \hfill $\Box$ 
  
\begin{coro} (a) If $z, w \in \mathds{Z}[\theta]$ and $(z) = (w)$, then, $\eta^+(z) = \eta^+(w)$.   
(b) If $(z)$ is an ideal with $\eta^+(z) > 0$, then, $z$ should be norm minimum among the elements of the ideal of non null norm.  
(c) If $w \in  (z)$ and $\eta^+(w) = \eta^+(z) > 0$, then, $(z) = (w)$, that is, $z$ and $w$ are associated elements.   
\end{coro}

\noindent{\bf Proof:}
\begin{description}
\item[a] Since $z = uw$ with $u$ unit, then $\eta^+(z) = \eta^+(uw) = \eta^+(u) \eta^+(w) = \eta^+(w)$ because $\eta^+(u) = 1$.  
\item[b] If $w \in (z)$ with $\eta^+(w) > 0$, then, $w = rz$, therefore, $\eta^+(w) = \eta^+(rz) = \eta^+(r)\eta^+(z) \geq \eta^+(z)$ because $\eta^+(r) \geq 1$ since it is non-zero $(\neq 0)$.  
\item[c] Since $w \in (z)$ then $w = rz$, therefore, $\eta^+(w) = \eta^+(rz) = \eta^+(r) \eta^+(z) = \eta^+(z)$ and, since $\eta^+(z) > 0$, we should have $\eta^+(r) = 1$, that is, $r$ is unit. \hfill $\Box$ 
\end{description}

\section{Primes Elements, Irreducible Elements and Factorization in $\mathds{Z}[\theta]$}  
In this section we adopt the usual definitions of \emph{prime element} and \emph{irreducible element} in a ring. An element $p \in \mathds{Z}[\theta]$  is called \emph{prime} if $p$ is non-zero, non unit and $p|xy$ implies $p|x$ or $p|y$. An element $a \in \mathds{Z}[\theta]$ is called \emph{irreducible} if $a$ is non-zero, non unit and $a = zw$ implies that $z$ is unit or $w$ is unit. It can be  proved that if $\mathcal{A}$ is an integral domain, then, all the prime elements of $\mathcal{A}$ are irreducible, and it can also be proved that if $\mathcal{A}$ is an unique factorization domain, in particular if it is Euclidean, then, all the irreducible elements are prime. In fact, the domain of Gaussian integers $\mathds{Z}[i]$ is Euclidean, consequently, the prime elements coincide with irreducible elements in it. In this paper, we study the \emph{prime} and \emph{irreducible} elements of $\mathds{Z}[j]$ and $\mathds{Z}[k]$, distinguishing the cases where they are zero divisors from the cases where they are not.   

In $\mathds{Z}[\theta]$ we have the following result:  
  
\begin{prop} 
If $a \in \mathds{Z}[\theta]$ and $a \notin \mathcal{D}$, then, $a$ prime implies $a$ irreducible.  \label{prop 5}
\end{prop}
  
\noindent{\bf Proof:} Suppose that $a$ is prime and $a = xy$, then, $a|xy$; therefore, $a|x$ or $a|y$, that is, $x = ra$ or $y = sa $ for some $r, s \in \mathds{Z}[\theta]$.

If $x = ra$, then, $a = xy = ray$, therefore, $(1 - ry)a = 0$. Since $a \notin \mathcal{D}$ we have that $ry = 1$, so, $y$ is unit. Similarly, if $y = as.$ \hfill $\Box$ 

Next we show that in $\mathds{Z}[j]\backslash \mathcal{D}$ and $\mathds{Z}[k] \backslash \mathcal{D}$ there exist irreducible elements that are not primes, therefore, in those cases, the set of prime elements in $\mathds{Z}[\theta] \backslash \mathcal{D}$ is a proper subset of the set of irreducible elements. In both cases the examples are the same: let us consider $c = 2$; we have that $\eta^+(c) = 4 > 0$. On the other hand, $c|0$ and $0 = (1 + j)(1 - j)$ in $\mathds{Z}[j]$ and $0 = k^2$ in $\mathds{Z}[k]$. However, it is easily verified that $c$ does not divide $1 + j, 1 - j$ and $k$. For example, if $2|1 + j$, then, $1 + j = 2(x + jy)$, so $2x = 1$ resulting $x =\dfrac{1}{2} \notin \mathds{Z}$. Therefore $2$ is not prime. In the same way it can be proved that no prime integer $p$ is prime in $\mathds{Z}[\theta]$ for $\theta = j$ or $k$, contrasting to $\mathds{Z}[i]$, where the prime integers $p$ such  that $p \equiv 3(\hspace{-8pt}\mod 4)$ are prime elements. However, we are going to see, further, that no odd prime of $\mathds{Z}$ is irreducible in $\mathds{Z}[j]$, although, every prime element of $\mathds{Z}$ is irreducible in $\mathds{Z}[k]$. 

On the other hand, we see that, in the studied rings, $2$ is irreducible. Suppose that $2 = ab$. Then, $\eta^+(a)\eta^+(b) = \eta^+(ab) = \eta^+(2) = 4$; therefore, $(\eta^+(a), \eta^+(b)) = (1, 4), (4, 1)$ or $(2, 2)$. In the first two cases $a$ or $b$ are unit. Now we have to prove that the third case is impossible. Suppose that $a = x + \theta y$ with $\theta = j$ or $k$ and $\eta^+(a) = 2$, which means that $|x^2 - \theta^2y^2| = 2$. In the case $\theta = k$ we have $x^2 = 2$; that is impossible for $x\in \mathds{Z}$. In the case $\theta = j$ we have $|x^2 - y^2| = 2$. In that case, $x$ and $y$ are even or $x$ and $y$ are odd. In any case, it is entailed that $4$ divides $|x^2 - y^2|$. That it is impossible, since $4$ does not divide $2$. Contrasting the exposed, it is verified easily that $2$ is reducible in $\mathds{Z}[i]$, since $2 = (1 + i)(1 - i)$, i.e, none of the factor is unit.  
  
\begin{teor}[Factorization Theorem in Product of Irreducible Elements] 
If $a \in \mathds{Z}[\theta]$ and $a \notin \mathcal{D} \cup U(\mathds{Z}[\theta])$, then,  there exist $u \in U(\mathds{Z}[\theta])$ and irreducible $q_1 , ... , q_m$ such that $a = uq_1 ... q_m$.
\end{teor}

\noindent{\bf Proof:} Let $C = \{\eta^+(z) | z \notin \mathcal{D} \cup U(\mathds{Z}[\theta])\}$. Then, $\eta^+(a) \in  C$ and for all $\eta^+(z) \in C, \eta^+(z) > 1$. The test will be made by induction on $\eta^+(a)$.\\
Step Base: $\eta^+(a) = {\rm min} C$. We will prove that $a$ is irreducible.  

Suppose that $a$ is not irreducible, then, $a = zw$ with $z$ and $w$ not unit.

Then, we have $\eta^+(a) = \eta^+(zw) = \eta^+(z)\eta^+(w)$ with $\eta^+(z) > 1$ and $\eta^+(w) > 1$, that is, $z, w \in C$. But, $\eta^+(z) < \eta^+(z)\eta^+(w) = \eta^+(a)$ and also $\eta^+(w) < \eta^+(a)$. That is an absurd by means of the minimality of $\eta^+(a)$.

Inductive step: $\eta^+(a) > {\rm min} C$.   

If $a$ is irreducible, there is nothing to demonstrate.  

Suppose that $a$ is not irreducible; then, $a = zw$ with $z$ and $w$ not unit. In fact, as in the argument above, $1 < \eta^+(z) < \eta^+(a)$ and $1 < \eta^+(w) < \eta^+(a)$. Therefore, by inductive hypotheses, $z = up_1 ... p_r$ and $w = vq_1 ... q_s$, with $u$ and $v$ unit and $p_1, ..., p_r$ and $q_1 , ..., q_s$  irreducible; therefore, $a = uvp_1 ... p_rq_1 ... q_s$ and also $uv$ is unit. \hfill $\Box$ 
  
It can easily shown that if the factorization of an element of $\mathds{Z}[\theta]$ is expressed by means of the primes, then, the factorization is unique. Therefore, in $\mathds{Z}[i]$ there exist a unique factorization. Regarding the non uniqueness of the factorization in irreducible elements in the rings $\mathds{Z}[j]$ and $\mathds{Z}[k]$ we can use the following examples. In $\mathds{Z}[j]$, consider $a = 8$. From the previous analysis, $8 = 2^3$ is a factorization into irreducible elements in $\mathds{Z}[\theta]$ for $\theta = j$ or $k$. We realize that $8 = (3 + j)(3 - j)$ is also a factorization into irreducible elements in $\mathds{Z}[j]$. If $3 \pm j$ = $zw$ with $z$ and $w$ not unit, then, $\eta^+(z) \eta^+(w) = \eta^+(zw) = \eta^+(3 \pm j) = 8$. So $\eta^+(z) = 2$ and $\eta^+(w) = 4$ or the opposite. But, as previously explained, $z \in \mathds{Z}[\theta]$ does not exist such that $\eta^+(z) = 2$. Therefore, $3\pm j$ is irreducible in $\mathds{Z}[j]$. In the case of $\mathds{Z}[k]$, a simple example of easy verification is the following: $4 = 2^2 = (2 + k)(2 - k)$. In fact, for similar considerations, $2 + k$ and $2 - k$ are irreducible elements in $\mathds{Z}[k]$.

\section{Characterization of the Prime Elements of $\mathds{Z}[\theta]$ } 
\begin{prop}
(a) If $p$ is prime in $\mathds{Z}[\theta]$, then, for all $z \in \mathcal{D}, p|z$ or $p|\bar z$. In particular, $p|1 + j$ or $p|1 - j$ in $\mathds{Z}[j]$ and $p|k$ in $\mathds{Z}[k]$.  
(b) In $\mathds{Z}[j]$, the elements $1 + j$ and $1 - j$ are prime, and in $\mathds{Z}[k], k$ is prime.  
\end{prop}

\noindent{\bf Proof:}
(a) Observe that if $z \in D$, then, $0 = \eta(z) = z\bar z$ and $p|0$.\\
(b) We have just going to prove that $1 + j$ is prime in $Z[j]$. Suppose $1 + j|ab$; then, $ab = z(1 + j)$ for some $z \in \mathds{Z}[j]$; therefore, $\eta^+(a)\eta^+(b) = \eta^+(ab) = \eta^+(z) \eta^+(1 + j) = 0$; therefore, $\eta^+(a) = 0$ or $\eta^+(b) = 0$. Suppose that $\eta^+(a) = 0$. Then $a \in \mathcal{D} = \mathcal{D}^+ \cup \mathcal{D}^-$, that is, $a \in \mathcal{D}^+$ or $a\in \mathcal{D}^-$. If $a \in \mathcal{D}^+$, then obviously, $1 + j |a$. If $a \in \mathcal{D}^-$, since $ab = 0$, we have that $b \in \mathcal{D}^+$; in that case, $1 + j |b.$\hfill $\Box$ 
  
The previous proposition shows that the prime elements of $\mathds{Z}[j]$  are divisors of $1 + j$ or $1 - j$, and the prime elements of $\mathds{Z}[k]$ are divisors of $k$. Next, we are about to see which elements are the divisors of $1 + j$ and $1 - j $ in $\mathds{Z}[j]$ and $k$ in $\mathds{Z}[k]$.  
  
\begin{prop} 
(a) $k$ is irreducible, therefore, the divisors of $k$ in $\mathds{Z}[k]$ are $k$ and its associates, that is, $\pm k$. On the other hand, $xk$ is reducible for all $x \in \mathds{Z}$ with $x \neq \pm 1$.  
(b) $1 + j$ and $1 - j $ are reducible and, therefore, all the non-zero elements of $\mathcal{D}$ in $\mathds{Z}[j]$ are reducible. On the other hand, the only divisors, up to the associated of $1 + j$ and $1 - j $ in $\mathds{Z}[j]$, but themselves, are the elements $(x + 1) - jx$, for all $x \in \mathds{Z}$, and $(x + 1) + jx$ for all $x \in \mathds{Z}$, respectively, where none of them is unit.  
\end{prop}

\noindent {\bf Proof:}\\
 \indent (a)	Suppose that $k = (x + ky)(x' + ky') = xx' + k(xy ' + x'y)$; therefore, $xx' = 0$ and $xy' + x'y = 1$.
  
If $x = 0$, then, $x'y = 1$; therefore, $x' = \pm 1$ for any $y'$; then, $x' + ky' = \pm 1 + ky' \in U(\mathds{Z}[k])$.
  
If $x\neq 0$, then, $x' = 0$; therefore, $xy' = 1$; therefore, $x = \pm 1$ for any $y$. Therefore, $x + ky = \pm 1 + ky \in U(\mathds{Z}[k])$.
  
On the other hand, if $x \neq \pm 1$, we cannot factorize $xk = (x + ky)k$, where none of the factors is unit, because $\eta(x + ky) = x^2 \neq 1$ and $\eta(k) = 0 \neq 1$.  

(b) We are going to prove the case $1 + j$: suppose that $1 + j = (x + jy)(r + js) = (xr + ys) + j(xs + yr)$; then, $xr + ys = 1$ and $xs + yr = 1$. Subtracting the previous equations, we have $(x - y)(r - s) = 0$; therefore either, $x = y$ or $r = s$. If $x = y$, then, $xr + xs = 1$, that is, $x(r + s) = 1$; therefore $x = \pm 1$ and $r + s = \pm 1$; therefore, $x = y = \pm 1$ and $r = \pm 1 - s$, so, if $x + jy = 1 + j$, then, $r + js = (1 - s) + js = (-s + 1) - j(-s)$, and if $x + jy = -1 - j$, then, $r + js = (-1 - s) + js = - ((s + 1) - js) .$ \hfill $\Box$ 
    
Due to the previous proposition, we have that the irreducible elements of $\mathds{Z}[k]$ contained in $\mathcal{D}$ are just $\pm k$, and that $\mathds{Z}[j]$ does not have irreducible elements in $\mathcal{D}$.  
  
\begin{coro}
The set of prime elements of $\mathds{Z}[k]$  is $\{z \in \mathds{Z}[k] | z = uk$ with $u$ unit$\} = \{\pm k\}$, that is, it coincides with the irreducible set of $\mathds{Z}[k]$ contained in $\mathcal{D}$.  
\end{coro}
  
\noindent{\bf Proof:}\\ 
$\subseteq$: If $p$ is prime, then $p$ is divisor of $k$; therefore, from the previous proposition, $p$ is associate of $k$.  \\
$\supseteq$: Since $k$ is prime, all the associates of $k$ are prime.  \hfill $\Box$ 

\begin{prop}
All the prime elements $p$ of $\mathds{Z}[\theta]\backslash\mathcal{D}$ divide one of the prime elements of the decomposition  of $\eta(p)$ in $\mathds{Z}$.  
\end{prop}

\noindent {\bf Proof:}
Let $p$ be a prime element of $\mathds{Z}[\theta]$ and let us consider the decomposition $p\bar p = \eta(p) = \pm 
q_1^{r_1} ... q_n^{r_n}$ with each $q_i$ prime; then, since $p|\eta(p)$, there exists $i$ such that $p|q_i$. \hfill $\Box$ 

\begin{prop}
If $\eta^+(a)$ is prime of $\mathds{Z}$, in particular, $a \notin \mathcal{D}$, then, $a$ is irreducible in $\mathds{Z}[\theta]$.\label{prop 9}
\end{prop}

\noindent{\bf Proof:} Suppose that $a$ is reducible, that is, $a = zw$ with $z$ and $w$ not units; then, $\eta^+(z) > 1$ and $\eta^+(w) > 1$; therefore, $\eta^+(a) = \eta^+(z) \eta^+(w)$ is not prime. \hfill $\Box$

The reciprocal of the previous proposition is not valid due to the elements $z = 2$ in the rings $\mathds{Z}[j]$ and $\mathds{Z}[k]$ and $z = 3$ in $\mathds{Z}[i]$.

We should observe that if $z = x + ky \in \mathds{Z}[k]$, then, $\eta^+(z) (=x^2)$ is never a prime. Therefore, the proposition above is not an approach to test if an element is irreducible in $\mathds{Z}[k]$. Besides, the mentioned fact can be used to prove that all the prime elements of $\mathds{Z}$ are irreducible in $\mathds{Z}[k]$. Let us suppose that $p$ is a reducible prime; then, there are $a, b$ not units in $\mathds{Z}[k]$ such that $p = ab$; therefore, $p^2 = \eta^+(p) = \eta^+(a)\eta^+(b)$, therefore, $\eta^+(a) = p (= \eta^+(b))$, which is impossible. On the other hand, due to the previous proposition, if $p$ is an odd prime of $\mathds{Z}$ of the form $2n + 1$, then, the element $z = (n + 1) \pm jn$ is irreducible in $\mathds{Z}[j]$, since $\eta^+(z) = |(n + 1)^2 - n^2| = |2n + 1| = |p|$.

The following theorem is essential in the characterization of the irreducible elements of $\mathds{Z}[i]$. It can also be proved for $\mathds{Z}[j]$ by replacing ``the sum of squares" by ``the difference of squares". However, that is unnecessary because, in fact, it is very frequent to find an integer number that is the difference of two squares. Indeed, every odd integer $2n + 1 = (n + 1)^2 - n^2$ is the difference between two squares; in particular all prime $p \neq 2$ is the difference between two squares, therefore, it is reducible in $\mathds{Z}[j]$. In addition, $p$ can be factored in the following way: $p = 2n + 1 = (n + 1)^2 - n^2 = ((n + 1) + jn)((n + 1) - jn)$, where both factors are irreducible and not unit. In contrast to that, $2$, as we have said before, is irreducible in $\mathds{Z}[j]$.   

\begin{prop}
Let $p$ be a prime of $\mathds{Z}$. The following statements are equivalent: \\
(i) $p$ is reducible in $\mathds{Z}[i]$.\\  
(ii) $p = a\bar a$, with $a$ irreducible in $\mathds{Z}[i]$.  \\
(iii) $p$ is the sum of two squares. 
\end{prop}
  
The following proposition characterizes the irreducibility of the elements of the form $(x + 1) \pm jx$, with $x\in \mathds{Z}$, which, as we have already expressed, are the divisors of $1 \pm j$ in $\mathds{Z}[j]\backslash\mathcal{D}$.  
  
\begin{prop}
Let $z = (x + 1) \pm jx$ with $x$ integer. The following statements are equivalent:\\
(i) $\eta(z)$  is prime in $\mathds{Z}$.\\  
(ii) $z$ is irreducible in $\mathds{Z}[j]$.\\  
(iii) $z$ is prime in $\mathds{Z}[j]$.  \label{prop 11}
\end{prop}

\noindent {\bf Proof:}\\
($i \Longrightarrow  ii$): It is an immediate consequence from Proposition \ref{prop 9}.  \\
($ii \Longrightarrow i$) Let us suppose that $\eta(z) = (x + 1)^2 - x^2 = 2x + 1 = ab$ with $a \neq \pm 1$ e $b \neq \pm 1$. Since $ab$ is odd, then, $a$ and $b$ are both odd, that is, $a = 2u + 1$ and $b = 2v + 1$.
  
Let us consider $z_1 = (u + 1) \pm  ju$ and $z_2 = (v + 1) \pm jv$. Therefore, $\eta(z_1) = 2u + 1 = a\neq\pm 1$ and $\eta(z_2) = 2v + 1 = b \neq\pm 1$, where, $\eta^+(z_1) \neq 1$ and $\eta^+(z_2) \neq 1$, that is, $z_1$ and $z_2$ are not unit. However, $z_1z_2 = ((u + 1)\pm ju)((v + 1) \pm jv) = (2uv + u + v + 1) \pm j(2uv + u + v)$. Furthermore, $2x + 1 = ab = (2u + 1)(2v + 1) = 4uv + 2u + 2v + 1$; therefore, $x = 2uv + u + v$ from which, $z_1z_2 = (x + 1) \pm jx = z$, a contradiction.  \\
($ii \Longrightarrow iii$) Let us suppose $z$ irreducible and $z|ab$. We will prove that $z|a$ or $z|b$. By hypothesis, there is $c$ such that $ab = cz $ ... (1). Let us suppose $a = a_1 + ja_2, b = b_1 + jb_2$ and $c = c_1 + jc_2$. Developing the identity (1) we obtain $a_1b_1 + a_2b_2 = c_1(x + 1) \pm c_2x$ and $a_1b_2 + a_2b_1 = c_2(x + 1) \pm c_1x$. Adding in the $+$ case and subtracting in the $-$ case we obtain, respectively, $(a_1 + a_2)(b_1 + b_2) = (c_1 + c_2)(2x + 1)$ and $(a_1 - a_2)(b_1 - b_2) = (c_1 - c_2)(2x + 1)$.

Since (ii) entails (i), we have that $2x + 1$ is prime in $\mathds{Z}$; therefore, in the $+$ case of: either $ 2x + 1|a_1 + a_2$ or $2x + 1|b_1 + b_2$, and in the $-$ case: either $2x + 1|a_1 -a_2$ or $2x + 1|b_1 - b_2$. Let us suppose that $2x + 1|a_1 + a_2$. We will prove that $z (= (x + 1) + jx)|a$. In fact, there is $r$ such that $a_1 + a_2 = r(2x + 1)$. Therefore, supposing $a_1 - a_2 = s$, we get $a_1 = rx + \dfrac{1}{2}(r + s)$ and $a_2 = rx + \dfrac{1}{2}(r - s)$. We may observe that $\dfrac{1}{2}(r + s) = a_1- rx$ and $\dfrac{1}{2}(r - s) = a_2 - rx$. Therefore, they are integer; where, $a = a_1 + ja_2 = (rx + \dfrac{1}{2}(r + s)) + j(rx + \dfrac{1}{2}(r - s)) = (\dfrac{1}{2}(r + s) + \dfrac{1}{2} j(r - s))((x + 1) + jx)$. So, $z|a$.\\
($iii \Longrightarrow ii$) It is straightforward from Proposition \ref{prop 5}, because $\eta(z) = 2x + 1 \neq 0$, that is, $z \notin \mathcal{D}.$ \hfill $\Box$ 

\begin{coro}
The set of prime elements of $\mathds{Z}[j]$ contained in $\mathcal{D}$ is $\{u(1 \pm j) | u$ is unit$\}$, while the set of prime elements of $\mathds{Z}[j]\backslash\mathcal{D}$ is $\{u((x + 1) \pm jx) | u$ is unit and $2x + 1$ is prime$\}.$  
\end{coro}

\section{Characterization of the Irreducible Non \\ Prime Elements of $\mathds{Z}[j]\backslash\mathcal{D}$}   
We saw, in the previous section, that $\mathcal{D}$ does not contain any irreducible of $\mathds{Z}[j]$. Let $\alpha \in \mathds{Z}[j]\backslash\mathcal{D}$ be a non prime irreducible, then, since $\eta^+(\alpha) \neq 0$, we have that $\eta^+(\alpha) = 2^\gamma\cdot p_1^{\gamma_1}\cdot \dots \cdot p_m^{\gamma_m}$ with the $p_k$ odd prime of $\mathds{Z}$.  
  
\begin{prop}
If $a$ is like above, then, $\gamma_1 = ... = \gamma_m = 0$; therefore, $\eta^+(a) = 2^\gamma$ for some $\gamma \geq 1$.
\end{prop}  

\noindent{\bf Proof:}
Let us suppose that $\gamma_k \geq 1$, then, $p_k$ is an odd prime, namely, $p_k = 2n + 1$, and $p_k|\eta^+(a)$. Let us consider $a_k = (n + 1) + jn$, then, from Propositions \ref{prop 9} and \ref{prop 11}, $a_k$ is prime in $\mathds{Z}[j]$ because $\eta^+(a_k) = p_k$. On the other hand, $a_k |a_k\overline{a_k} = p_k|\eta^+(a) =\pm\eta(a) = \pm a\bar a$. Therefore, interchanging, $a_k$ by $\overline{a_k}$ if necessary, we have that $a_k|a$.  Next, $a = \beta a_k$ with $\eta^+(\beta)\neq 1$ because $a$ is not prime, and so $a$ reducible, a contradiction.\hfill $\Box$

The next step is to find all the elements $a$ which are irreducible in $\mathds{Z}[j]\backslash\mathcal{D}$ such that $a\bar a = \pm 2^\gamma$ with $\gamma \geq 2$  (we saw already that $a$ does not exist when $a\bar a = \pm 2$). It is worth to observe that all $a \notin \mathcal{D}$ has an associated with $\eta(a) > 0$ and $Re(a) > 0$; hence we can suppose that fact. Let, then, $a$ be such that $a\bar a = 2^\gamma$ with $\gamma \geq 2$. Observing that $2^\gamma = (2^{\gamma -2} + 1)^2 - (2^{\gamma -2} - 1)^2$, we can take $a = (2^{\gamma -2} + 1) \pm  j(2^{\gamma -2} - 1)$. We will see that such $a$, up to associated elements, is the only irreducible which is non prime of $\mathds{Z}[j]$. For this, we need the following lemma of immediate verification.  

\begin{lema}
If $z$ and $w$ are of the from $(2n + 1) \pm j(2n - 1)$, then, $2|zw.$ 
\end{lema}

\begin{prop}
Let $a$ be such that $\eta(a) = 2^{\gamma+2}$ with $\gamma \geq 0$ e $Re(a) > 0$. Then, $a$ is irreducible iff (*) $a = (2^{\gamma} + 1) \pm  j(2^{\gamma} - 1)$.
\end{prop}   
  
\noindent{\bf Proof:}  \\
($\Longrightarrow $) Let us suppose that $a = x + jy$ is not in the form (*). We will prove that $2|a$. Since $\eta(a) = x^2 - y^2 = (x + y)(x - y) = 2^{\gamma+2}$, we have that $x + y = 2^{\gamma+2-h}$ and $x - y = 2^h$ with $0 \leq h \leq \gamma + 2$. Realize that, if $h = 0$ or $h = \gamma + 2$ we would have $x + y$ even (resp. odd) and $x - y$ odd  (resp. even). If $h = 1$ or $h = \gamma + 1$ we would have $a$ in the form (*). Therefore, $2 \leq h \leq \gamma$. Solving the system we have $x = 2^{\gamma +1-h} + 2^{h-1}$ and $y = 2^{\gamma+1-h} - 2^{h-1}$, from which $a = (2^{\gamma +1-h} + 2^{h-1}) + j(2^{\gamma +1-h} - 2^{h-1}) = 2[(2^{\gamma -h} + 2^{h-2}) + j(2^{\gamma -h} + 2^{h-2})]$, that is, $2|a$. Therefore, $a$ is reducible.  \\
($\Longleftarrow$) Let $a$ be in the form (*) and suppose that it is reducible, that is, $a = zw$ with $\eta^+(z) \geq 2$ and $\eta^+(w) \geq 2$. In fact, we can suppose $\eta(z) = 2^{\gamma+2-h}$ and $\eta(w) = 2^h$ with $1 \leq h \leq \gamma + 1$. Moreover, either $h = 1$ or $h = \gamma + 1$ are impossible as we already saw. Therefore, $2 \leq h \leq \gamma$. Since 2 does not divide $a$ due to its form, we have, from the lemma, that either $z$ or w are not in the form (*). Therefore, from $(\Longrightarrow)$ we have that $2|z$ or $2|w$, from which $2|zw (= a)$, a contradiction.\hfill $\Box$ 

\begin{coro}  
The irreducible not prime elements of $\mathds{Z}[j]\backslash\mathcal{D}$ are the associated of $(2n + 1) \pm j(2n - 1)$ for all $n \geq 0.$
\end{coro}
  
Finally, in $\mathds{Z}[j]$ we can enunciate and demonstrate the analogous of the \emph{ Fermat Theorem} that characterizes, in $\mathds{Z}[i]$, the positive integers that are sum of two squares.  
  
\begin{teor}
Let $n > 0$ be an integer with decomposition in primes factors given for  
$$n = 2^\gamma\cdot p_1^{\gamma_1}\cdot ... \cdot p_m^{\gamma_m},$$
then, $n$ is difference of two squares if, and only if, $\gamma \neq 1$.  
\end{teor}

\noindent{\bf Proof:}\\  
$(\Longrightarrow)$ If $\gamma = 1$, then, $n = 2p_1^{\gamma_1}\cdot ... \cdot p_m^{\gamma_m}$, and hence, $2|n = r^2 - s^2 = (r - s)(r + s)$, from which, either $2|(r - s)$ or $2|(r + s)$. In fact, $2$ divides both. Therefore, $4 (=2^2)|n$, a contradiction.  \\
$(\Longleftarrow)$ Let us suppose $\gamma \neq 1$. If $\gamma = 0$, then, since any odd prime is a difference of squares and this property is preserved by products, we have that $n$ is difference of squares. If $\gamma \geq 2$, then, since $2^\gamma = (2^{\gamma -2} + 1)^2 – (2^{\gamma -2} - 1)^2$, we have that, also, n is difference of squares.   \hfill $\Box$ 
  
\section{Characterization of the Irreducible Not Prime Elements of $\mathds{Z}[k]\backslash\mathcal{D}$}  
In the next proposition, we will make extensive use of the following elementary property about \emph{Diophantine equations}: if $a, b$ and $c$ are integer numbers, then, the equation $ax + by = c$ has integer solution if, and only if, ${\rm gcd} (a,b)|c$.  
 
\begin{prop}
Let $z \in \mathds{Z}[k]\backslash\mathcal{D}$, that is, $Re(z) \neq 0$ (we can suppose $Re(z) > 0$). Then:  \\
(a) If $z = p + ky$, with $p$ prime, then $z$ is irreducible.\\  
(b) If $z = x + ky$, with $x$ non prime's potency, then $z$ is reducible.  \\
(c) If $z = pg + ky$, with $p$ prime and $\gamma \geq 2$, then $z$ is reducible $\Longleftrightarrow p|y$.  
\end{prop}

\noindent{\bf Proof:}\\  
(a) Let $z = p + ky$. Let us suppose $z = ab$. Then, $p^2 = \eta(z) = \eta(a)\eta(b)$. Therefore, since $\eta(a)$ and $\eta(b)$ cannot be equal to $p$, the only possibility that we may have is either $\eta(a) = 1$ (and $\eta(b) = p^2)$ or $\eta(b) = 1$ (and $\eta(a) = p^2$), that is, either $a$ is unit or $b$ is unit. Therefore, $z$ is irreducible.  \\
(b) Let us suppose that $x$ is not a power of a prime. Then, $x = mn$ with $m \neq\pm 1, n \neq\pm 1$ and ${\rm gcd} (m,n) = 1$. Therefore, we may decompose $z = (m + kr)(n + ks) = x + k(rn + sm)$ and the equation $rn + sm = y$ has solution because ${\rm gcd} (m,n) = 1$.\\  
(c) $(\Longleftarrow)$ Let us suppose that $z = p^\gamma  + ky$, with $p$ prime, $\gamma \geq 2$ and $p$ divides $y$. Then, we can decompose $z = (p + kr) (p^{\gamma -1} + ks) = x + k(p^{\gamma - 1}r + ps)$ and the equation $p^{\gamma - 1}r + ps = y$ has solution because ${\rm gcd} (p^{\gamma - 1}, p) = p$ and $p|y$.\\  
$(\Longrightarrow)$ Let us suppose that $z = p^\gamma  + ky$, with $p$ prime, $\gamma \geq 2$ and $p$ does not divide $y$, and let us suppose $z = (p^{\gamma - h}  + kr)(p^h  + ks)$ with $0 \leq h \leq \gamma$. We will prove that either $h = 0$ or $h = \gamma$, in which case some of the factors are unit and, therefore, $z$ would be irreducible. In fact, if $1\leq h \leq \gamma -1$, we would have ${\rm gcd} (p^{\gamma - h} , p^h) = p$, therefore, since $p$ does not divide $y$, the equation $p^{\gamma - h}s + p^hr = y$ would have no solution, a contradiction. \hfill $\Box$

\end{document}